\documentclass{article}
\usepackage{amssymb, amsmath, amsfonts}
\usepackage{bm}
\usepackage[all,2cell]{xy}\UseAllTwocells

\def\genfd{{\bf k}} 
\def\id{{\rm id}}
\def\Id{{\rm Id}}
\def\AA{{\mathcal A}}
\def\BB{{\mathcal B}}
\def\CC{{\mathcal C}}
\def\DD{{\mathcal D}}

\def\MM{{\mathcal M}}

\def\tto{\Rightarrow}

\def\nxpoint{\refstepcounter{subsection}%
  \makepoint{\thesubsection}}
\def\nxsubpoint{\refstepcounter{subsubsection}%
  \makepoint{\thesubsubsection}}
\def\refpoint#1{{\rm\textbf{\ref{#1}}}}

\def\makepoint#1{\medbreak\noindent{\bf #1. }}
\newcommand{\nodo}[1]{}

\begin{document}
\begin{center}
{\bf Compatibility of (co)actions and localizations}

{\sc Zoran \v{S}koda}, {\tt zskoda@irb.hr}

{\it preliminary version}
\end{center}

{\bf Earlier, Lunts and Rosenberg studied a notion of
compatibility of endofunctors with localization functors,
with an application to the study of differential operators on
noncommutative rings and schemes. Another compatibility
-- of Ore localizations of an algebra with a comodule algebra
structure over a given bialgebra -- introduced in my earlier work
-- is here described also in categorical language,
but the appropriate notion differs from that of Lunts and Rosenberg,
and it involves a specific kind of distributive laws. Some basic facts about
compatible localization follow from more general functoriality 
properties of associating comonads or even actions of monoidal categories
to comodule algebras. We also introduce localization 
compatible pairs of entwining structures.
}

\vskip .2in

\section{Introduction}

\nxpoint {\bf (Notation and prerequisites).}
Throughout the paper $\genfd$ is a ground field, but for most
results it can be taken to be a commutative ring.
The unadorned tensor symbol means tensoring over $\genfd$.
We assume that the reader is well familiar
with adjoint functors and familiar with (co)monads
(\cite{Borceux,MacLane:CWM,skoda:nloc}) which some call triples
(\cite{Semtriples}). We will speak (co)modules over (co)monads
for what many call (co)algebras over (co)monads.

\nxpoint {\bf (Context and motivation).} 
Apart from general purpose, this article
is aimed to create the preliminaries for a natural general theorem
on reconstruction of the structure of the {\it noncommutative
scheme} (\cite{Ros:NCSch}) on the category of equivariant
quasicoherent sheaves over a noncommutative scheme which is
locally (in the sense of a cover by biflat affine localizations
compatible with (co)actions) of Galois type (noncommutative
principal bundle). This theorem involves the construction and
exactness properties of various adjoint pairs of
functors~(\cite{Ros:NCSch}), what asked for precise requirements
and usage of the correct and natural compatibility properties of
localization functors with (co)actions in such geometric setup.
This is the main subject of our forthcoming
paper~\cite{skoda:globGalois}.

\nxpoint {\bf (Localization.)} (\cite{GZ,skoda:nloc,Ros:NCSch}) A
{\bf localization functor} is a functor which is universal among
the functors inverting a given class of morphisms in a domain
category. A {\bf continuous localization} functor is a functor
$Q^*:\AA\to\BB$ having a fully faithful right adjoint
$Q_*:\BB\to\AA$ (this implies that $Q^*$ is a localization
functor). This is equivalent to having a pair of adjoint functors
$Q^*\dashv Q_*$ for which the counit $\epsilon:Q^* Q_*\to\id_\BB$
is an isomorphism of functors. Consequently the multiplication
$Q_*\epsilon Q^* : Q_*Q^* Q_*Q^*\to Q_* Q^*$ of the monad induced
by this adjunction is clearly also an isomorphism (``idempotent
monad'') and the localized category $\BB$ is via the comparison
functor $N\mapsto (Q_* N, Q_*(\epsilon))$ equivalent to the
(Eilenberg-Moore) category of modules over that monad. One usually
says that $\BB$ is a {\bf reflective subcategory} of $\AA$
(strictly full subcategory where the inclusion has a left
adjoint).

In Abelian categories, one usually considers (additive) flat (=
exact and continuous) localization functors: they may be obtained
by localization at {\bf thick subcategories} (= full subcategories
closed under direct sums, quotients and extensions). The main
example in the categories of modules is any {\it Gabriel
localization} (at a Gabriel filter) of the category of left
modules over a unital ring $R$. Even better subclass is the class
of {\bf Ore localizations}, which are of the form $Q^* M = S^{-1}
R \otimes_R M$ where $S^{-1}R$ is the Ore localized ring, at a
(say left) Ore set $S\subset R$. In that case, $Q^*$ and $Q_*$ are
exact, $S^{-1}R$ is consequently flat over $R$ and the component
of the unit of the adjunction $\eta : \id\to Q_* Q^*$ for the ring
$R$ itself, namely $\iota_S : R\to S^{-1}R$, is a morphism of
unital rings. The multiplication induces an isomorphism of
$S$-modules $S^{-1}R\otimes_R S^{-1}R\to S^{-1}R$ (because the
monad $M\mapsto S^{-1}R\otimes_R M = Q^* Q_* M$ is an idempotent
monad).

\section{Functoriality of actegories from comodule algebras}

\nxpoint Let $B$ be a $\genfd$-bialgebra. The category ${}_B\MM$ 
of left $B$-modules is a monoidal category in standard way: 
the tensor product is the tensor product of the 
underlying $\genfd$-modules with the left $B$-action given by 
$b(x\otimes_\genfd y) = \sum b_{(1)}x\otimes_\genfd b_{(2)}y$.

\nxpoint Given a $\genfd$-bialgebra $B$,
a left (right) $B$-{\bf comodule algebra}
is a pair $(A,\rho)$ of an algebra $A$
and a left (right) $B$-coaction $\rho : A\to B\otimes A$
(resp. $\rho : A\to A\otimes B$) which is an algebra map.
We use extended Sweedler notation
$\rho(e) = \sum e_{(0)}\otimes e_{(1)}$ (\cite{Majid}). 

\nxsubpoint ${}_B\MM$ acts on ${}_\genfd\MM$ in a trivial way: on objects 
just tensor the underlying $\genfd$-modules; bialgebra $B$ lives in
the category of $\genfd$-modules, and its meaning is related to the
tensor product in ${}_\genfd\MM$. Thus this distinguished
{\it defining} or {\bf base action} 
is however important, because in noncommutative geometry 
it is natural to consider actions of ${}_B\MM$ which are 
{\bf geometrically admissible}. 
These are the actions of the type 
$\triangleleft : \CC\times{}_B\MM\to \CC$ on an abstract category $\CC$
equipped with a direct image functor $U:\CC\to {}_\genfd\MM$ 
such that $U\circ \triangleleft=\triangleleft_0 \circ (U\times\Id_{{}_B\MM})$,
where $\triangleleft_0$ is the base action 
$\triangleleft_0:{}_\genfd\MM\times{}_B\MM\to{}_\genfd\MM$.
Such actions may be called {\it lifts} of $\triangleleft_0$ along $U$. 
Lifts to $\CC= {}_E\MM$ where $E$ is a $\genfd$-algebra 
are in a bijective correspondence with the distributive laws 
between the base action and monad $E\otimes_\genfd$ on ${}_\genfd\MM$
)\cite{skoda:distr,skoda:gmj}. Such distributive laws are generalizations of
Beck' classical distributive laws between two (co)monads. 

\nxsubpoint The distributive law with components 
$l_{M,Q} : E\otimes (M\otimes Q)\to (E\otimes M)\otimes Q$ 
given by $e\otimes (m\otimes q)\mapsto\sum e_{(0)} \otimes m\otimes e_{(1)} q$ 
where $e\in E, m\in M, q\in Q$, where $M\in {}_\genfd\MM$ 
and $Q\in {}_B\MM$ induces thus a right action of ${}_B\MM$ 
on ${}_E\MM$ lifting the base action. More explicitly,
for $N\in{}_E\MM$, $N\triangleleft Q$ is $N\otimes Q$ with
left $E$-module structure $e(n\otimes q) = \sum e_{(0)} n\otimes e_{(1)}q$.

We want to describe in this section the functoriality of this and various
dual constructions. We often call monoidal category $\DD$ together with a 
(left or right) $\DD$-action on some category $\CC$ 
a (left or right) $\DD$-actegory.

\nxpoint \label{s:Hopfmodules} {\bf (Comonad for the relative Hopf
modules).} $B$ is a comonoid in the monoidal category ${}_B\MM$. 
Therefore the strong monoidal action of ${}_B\MM$ 
on any category sends it to a comonoid in the category of
endofunctors (in our case also additive). 
The underlying endofunctor $G : {}_E\MM\to{}_E\MM$ in the category
${}_E\MM$ of left $E$-modules on objects $M$ in $\MM$ is given
by the formula $G: M\mapsto M\otimes B$, where
the left $E$-module structure on $M\otimes B$ is given by
$e(m\otimes b):= \rho_E(e)(m\otimes b)= \sum e_{(0)}
m\otimes e_{(1)}b$ where $e\in E, m\in M, b\in B$. The
comultiplication $\Delta_B$ on $B$ induces the comultiplication
$\delta =\id\otimes\Delta:G\to GG$ on $G$ with counit
$\epsilon^G=\id\otimes\epsilon$ making
$\mathbf{G}=(G,\delta,\epsilon^G)$ a comonad (cf. the coring picture
in~\cite{BrzWis:corings}).

It is well-known that the category $({}_E\MM)_{\bf G}$ 
of $\mathbf{G}$-comodules (coalgebras) is equivalent to 
the category ${}_E\MM^B$
of left-right relative $(E,B)$-Hopf modules; thus we say that 
${\bf G}$ is the {\it comonad for relative Hopf modules}.
Left-right {\bf relative $(E,B)$-Hopf module} is a m
left module $(N,\nu_N)$ where
$\nu_N:E\otimes N\to N$ is a left $E$-action, equipped with a
right $B$-coaction $\rho_N:N\to
N\otimes B$ such that 
$\rho_N(\nu(e,n))=(\nu\otimes\mu_B)(\id\otimes\tau_{B,N}\otimes\id)
(\rho_E(e)\otimes\rho_N(n))$ for all $e\in E$, $n\in N$;
here $\tau_{B,N}:B\otimes N\to N\otimes B$
is the flip of tensor factors. Morphisms of relative Hopf modules
are morphisms of underlying $\genfd$-modules, which
respect $E$-actions and $B$-coactions.

\nxpoint We will work in part of the article 
not only with maps of comodule algebras over a fixed bialgebra, 
but we will also allow variable bialgebras.
Thus let $\phi:B\to B'$ be a map of bialgebras, $(E,\rho)$ a $B$-comodule
algebra and $(E',\rho')$ a $B'$-comodule algebra. 
Then a map of underlying algebras $f: E\to E'$ 
is a {\bf map of comodule algebras over} $\phi$ if
$\rho'\circ f= (f\otimes\phi)\circ\rho : E\to E'\otimes B'$. Alternatively,
one says that the pair $(f,\phi)$ is a 
{\it map of comodule algebras over varying bialgebras}. 

\nxpoint \label{p:lifteCompExtScalars}
For any algebra map $f:A\to A'$ we use geometric 
{\it inverse image} notation for 
the extension of scalars $f^* : M\to f^*M = A'\otimes_A M$
for the categories of left modules, though this implies
that $f\mapsto f^*$ is a covariant functor. 

\nxpoint \label{s:basicThm}{\bf Theorem.} {\it There is a canonical 2-cell 
$$\xymatrix{
{}_E\MM\times{}_B\MM\ar[r]^{f^*\times\phi^*}
\ar[d]^{\triangleleft}
&
{}_{E'}\MM\times{}_{B'}\MM\ar[d]^{\triangleleft'}
\\
{}_E\MM\ar[r]^{f^*}
& {}_{E'}\MM\ultwocell<\omit>{<0>\alpha}
}$$
that is a natural transformation
$\alpha = \alpha^{f,\phi} : 
f^*\circ\triangleleft \tto \triangleleft'\circ (f^*\times\phi^*)$.
}

{\it Proof.} The components 
$$\alpha_{M,Q}=\alpha_{M,Q}^{f,\phi} :E'\otimes_E(M \triangleleft Q)\to 
(E'\otimes_E M)\triangleleft' (B'\otimes_B Q)
$$ 
of the natural transformation $\alpha$,
where $M\in {}_E\MM$ and $Q\in {}_B\MM$ are defined as 
$\genfd$-linear extensions of the formulas
$$
\alpha_{M,Q}(e'\otimes_E (m\otimes q)) = \sum (e'_{(0)}\otimes_E m)\otimes
(1_{B'}\otimes_B e'_{(1)} q).
$$
One checks that $\alpha_{M,Q}$ is well defined (it is well-defined
before we quotient to $\otimes_E$; consider values
on $e'e\otimes (m\otimes q)$ and $\sum e'\otimes e_{(0)}m\otimes e_{(1)}q$ 
and calculate that both give the same) and 
that $\alpha_{M,Q}$ is indeed a morphism in ${}_{E'}\MM$.

\nxpoint \label{s:comonadMap}
Consider now the comonad for Hopf modules 
${\bf G} = (G,\delta,\epsilon)$ on ${}_E\MM$.

{\bf Corollary.} {\it  There is a 2-cell
$$\xymatrix{
{}_E\MM\ar[r]^{f^*}
\ar[d]^{G_E}
&
{}_{E'}\MM\ar[d]^{\triangleleft'}
\\
{}_E\MM\ar[r]^{f^*}
& {}_{E'}\MM\ultwocell<\omit>{<0>\alpha_B}
}$$
which is in fact a morphism of comonads. 
}

\nxpoint {\bf Theorem.} {\it 2-cells $\alpha$ paste correctly with respect
to composition of comodule algebra maps over varying base.
In other words, for components at $B$, the pasting
$$\xymatrix{
{}_E\MM\ar[r]^{f^*}\ar@/^2pc/[rr]^{(gf)^*}
\ar[d]_{G_E}
&
{}_{E'}\MM\ar[d]_{G_{E'}}\ar[r]^{g^*}
&
{}_{E''}\MM\ar[d]^{G_{E''}}
\\
{}_E\MM\ar[r]_{f^*}\ar@/_2pc/[rr]^{(gf)^*}
& {}_{E'}\MM\ultwocell<\omit>{<0>\alpha_B^{f,\phi}}\ar[r]_{g^*}
& {}_{E''}\MM\ultwocell<\omit>{<0>\alpha_{B'}^{g,\phi'}}
}$$
equals the two cell
$$\xymatrix{
{}_E\MM\ar[rr]^{(gf)^*}
\ar[d]_{G_E}
&&
{}_{E''}\MM\ar[d]^{G_{E'}}
\\
{}_E\MM\ar[rr]_{(gf)^*}
&& {}_{E''}\MM\ultwocell<\omit>{<0>\alpha_B^{gf,\phi'\phi}}
}$$
where the symbols for
canonical invertible 2-cells $g^* f^* \cong (gf)^*$ are ommitted.

Similar diagrams hold for other components.
}

Proof is an easy direct calculation.

\nxpoint Now we use the transformation $\alpha$ to 
induce the map for the categories of Hopf modules
${}_E\MM^B\to {}_{E'}\MM^{B'}$. It is known
that ${}_E\MM^B\cong ({}_E\MM)_{{\bf G}_E}$ so this procedure is standard.

Let $M\in {}_E\MM$ and $\rho_M:M\to M\otimes B$ be a coaction making
$M$ a relative Hopf module. 

{\bf Proposition.} {\it The extension of scalars $f^*:{}_E\MM\to{}_{E'}\MM'$  
lifts to the functor $f^{*B}:{}_E\MM^B\to {}_{E'}\MM^{B'}$ between the
categories of relative Hopf modules,
which is at objects given by 
$f^{*B}: (M,\rho_M)\mapsto (f^* M,\alpha_M\circ f^*(\rho))$. 
}

\section{Compatibility for comodule algebras}

\nxpoint {\bf (Compatibility of coactions and localizations).}

Let $(E,\rho)$ be a right $B$-comodule algebra.
An Ore localization of rings $\iota_S : E\to S^{-1}E$
is $\rho$-{\bf compatible} (\cite{skoda:ban})
if there exist an (automatically unique)
coaction $\rho_S : S^{-1}E\to S^{-1}E\otimes B$ making
$S^{-1}E$ a $B$-comodule algebra, such that the diagram
$$\xymatrix{
E\ar[r]^\rho\ar[d]^{\iota_S}& E\otimes B\ar[d]^{\iota_S\otimes B}
\\
S^{-1}E \ar[r]^{\rho_S}&S^{-1}E\otimes B
}$$
commutes. It is easy to check that this is equivalent to
an effective criterium that
for all $s\in S$, $(\iota_S\otimes\id_B)\rho(s)$
is invertible in $S^{-1}E\otimes B$. If the Ore localization is
$\rho$-compatible, $\rho_S$ is called (the)
{\bf localized coaction}. The elements $u\in S^{-1}E$ satifying
$\rho_S(u)= u\otimes 1$ i.e. the coinvariants under the
localized coaction are called {\bf localized coinvariants}. It is a basic
and important observation that the localization and taking coinvariants
do not commute: the subalgebra $(S^{-1}E)^{\mathrm{co}B}\subset S^{-1}E$
of localized coinvariants typically contains some extra elements which
do not naturally belong to the $\genfd$-submodule
$\iota_S(E^{\mathrm{co}B})$; moreover typically $\iota_S$ restricted to
the subring $E^{\mathrm{co}B}\subset E$
is not underlying a ring localization
$U^{-1} E^{\mathrm{co}B}$ with respect to any Ore subset $U$
in $E^{\mathrm{co}B}$.

\nxpoint \label{s:compOreCateg}{\bf Theorem.} {\it Let $B$ be a
$\genfd$-bialgebra, $(E,\rho)$ a $B$-comodule algebra,
$\mathbf{G}$ a comonad from \refpoint{s:Hopfmodules}, and $\iota :
E\to E_\mu$ a perfect (e.g. Ore) localization of rings, which
happens to be $\rho$-compatible. 

The $\genfd$-linear map
$$
l_{M,P} : E_\mu\otimes_E(M\otimes P)\tto (E_\mu\otimes_E M)\otimes P,
\,\,e\otimes(m\otimes b)\mapsto \sum(e_{(0)}\otimes m)\otimes
e_{(1)},
$$
for $m\in M, p\in P, e\in E$,
where $P$ is a $B$-module and $M$ a $E$-module is a well-defined morphism
of left $E$-modules. 
All $l_{M,P}$ together form a mixed distributive
law between the localization monad
$Q_* Q^*$ and the categorical action of ${}_B\MM$ on ${}_E\MM$. 
}

{\it Proof.} This is a slight generalization of the case $P=B$ 
which gives the distributive law between  the localization monad
$Q_* Q^*$ and the comonad $\mathbf G$ which is proved in~\cite{skoda:nloc}.
The general proof is analogous. 

\nxpoint {\bf Proposition.} {\it Given any continuous localization
functor $Q^*:\AA\to\AA_\mu$ and a comonad $\bf G$ together with
any mixed distributive law $l:Q_* Q^* G\tto G Q_* Q^*$,

1) $G_\mu = Q^* G Q_*$ underlies a comonad ${\bf
G}_\mu=(G_\mu,\delta^\mu, \epsilon^{G_\mu})$ in $\AA_\mu$ with
comultiplication $\delta^\mu$ given by the composition
$$
Q^* G Q_* \stackrel{Q^* \delta^G Q_*}\longrightarrow Q^* GG Q_*
\stackrel{Q^* G\eta G Q_*}\longrightarrow Q^* G Q_* Q^* G Q_*
$$
and whose counit $\epsilon^{G_\mu}$ is the composition
$$\xymatrix{
Q^*GQ_* \ar[r]^{Q^*\epsilon^G
Q_*}&Q^*Q_*\ar[r]^{\epsilon}&\Id_{\AA_\mu} }$$ (where the
right-hand side comultiplication $\epsilon$ is the counit of the
adjunction $Q_*\dashv Q^*$).

2) the composition
$$\xymatrix{
Q^* GM\ar[rr]^{Q^*(\eta_{GM})}&& Q^* Q_* Q^* GM
\ar[rr]^{Q^*(l_M)}&& G_\mu Q^* M }$$ defines a component of a
natural transformation $\alpha = \alpha_l: Q^* G\tto G_\mu Q^*$
for which the mixed pentagon diagram of transformations
$$\xymatrix{
Q^* G\ar[rr]^{\alpha}\ar[d]_{Q^*\delta^G}&& G_\mu Q^*\ar[d]^{\delta^\mu Q^*}\\
Q^* GG \ar[r]^{\alpha G}&G_\mu Q^* G\ar[r]^{G\alpha}& G_\mu G_\mu
Q^* }$$ commutes and $(\epsilon^{G_\mu} Q^*)\circ\alpha = Q^*
\epsilon^G$. In other words, $(Q^*,\alpha_l): (\AA, {\bf G})\to
(\AA_\mu, {\bf G}_\mu)$ is (up to orientation convention which
depends on an author) a map of comonads
(\cite{skoda:ent,Street:formal}). }

\nxpoint {\bf Theorem.} {\it Under assumptions in
\refpoint{s:compOreCateg}, there is a unique induced continuous
localization functor $Q^{B*}: {}_E\MM^B\to {}_{E_\mu}\MM^B$
between the categories of relative Hopf modules such that $U_\mu
Q^{B*} = Q^* U$ where $U$ and $U_\mu$ are the forgetful functors
from the category of relative Hopf modules to the categories of
usual modules over $E$ and $E_\mu$ respectively. }

{\it Proof.} We have stated this theorem and given direct proof in 
\cite{skoda:gmj}. The more general results from the previous sections make
it a special case of \refpoint{p:lifteCompExtScalars}.

\nxpoint {\bf Corollary.} {\it Let $E=B$ and $\genfd$ is a field. 
The only $\Delta$-compatible
Ore localization $B\to B_\mu$ is the trivial one.
}

This is an analogue of the statement that the only $G$-invariant Zariski open
subset of an algebraic group over a field is the whole group.

{\it Proof.} The compatible localization functor induces in this case
a localization functor $Q^{B*}: {}_B\MM^B\to {}_{B_\mu}\MM^B$. 
By the fundamental theorem on relative Hopf modules 
the domain of this functor is ${}_B\MM^B\cong {}_\genfd\MM$. 
But ${}_\genfd\MM$ is
just a category of vector spaces over a field which does not have nontrivial
continuous localizations with contradiction. 

\section{Compatibility for entwinings}

\nxpoint ({\bf Localization-compatible pairs of entwinings.}) Let
$A$ be a $\genfd$-algebra, $C$ a $\genfd$-coalgebra, $\iota: A\to
A_\mu$ a perfect localization of rings, and $\psi : A\otimes C\to
C\otimes A$, $\psi_\mu : A_\mu\otimes C\to C\otimes A_\mu$
entwinings. We say that $(\psi,\psi_\mu)$ is $\iota$-compatible
pair of entwinings if the diagram
$$\xymatrix{
A\otimes C\ar[r]^\psi\ar[d]^{\iota\otimes C}&C\otimes A\ar[d]^{C\otimes\iota}
\\
A_\mu\otimes C\ar[r]^{\psi_\mu}&C\otimes A_\mu
}$$
commutes.

\nxpoint Define comonad $\bf G$ on ${}_A\MM$ as usual: $G(M,\nu) =
(C\otimes M, (C\otimes\nu) \circ \psi_M)$.\vskip .04in

{\bf Proposition.} {\it Given a $\iota$-compatible pair
$(\psi,\psi_\mu)$ of entwinings, the $\genfd$-linear map
$\psi_\mu\otimes M : A_\mu \otimes C\otimes M\to C\otimes A_\mu
\otimes M$ factors to a well-defined map of $A$-modules}
$$
l_M : A_\mu\otimes_A GM \to G(A_\mu\otimes_A M).
$$
{\it Proof.} Consider the diagram
\begin{equation}\label{eq:doubleeq}\xymatrix{
A_\mu\otimes A\otimes C\otimes M \ar[d]
\ar@<-.6ex>[rr]_-{A_\mu\otimes (C\otimes\nu)\circ(\psi\otimes M)}
\ar@<.6ex>[rr]^-{\mu\otimes C\otimes M}&&A_\mu\otimes C\otimes
M\ar[r]\ar[d]^{\psi_\mu\otimes M}&A_\mu\otimes_A
GM\ar@{.>}[d]^{l_M}\\
C\otimes A_\mu\otimes A\otimes M \ar@<-.6ex>[rr]_-{A_\mu\otimes
(C\otimes\nu)\circ(\psi\otimes M)} \ar@<.6ex>[rr]^-{C\otimes
\mu\otimes M}&& C\otimes A_\mu\otimes M\ar[r]&G(A_\mu\otimes_A
M)}\end{equation} where the rows are equalizer forks and the left
vertical arrow is the composition
$$\xymatrix{
A_\mu\otimes E\otimes C\otimes
M\ar[r]\ar@<.6ex>@{}[r]^{A_\mu\otimes\psi\otimes M}&A_\mu\otimes
C\otimes E\otimes M\ar[r]\ar@<.6ex>@{}[r]^{\psi_\mu\otimes
A\otimes M}&C\otimes A_\mu\otimes A\otimes M}$$

If the left square in~(\ref{eq:doubleeq}) sequentially commutes,
then clearly the right vertical arrow factors to a well-defined
map $l$.

In the following two diagrams we omit the tensor product sign
$\otimes_\genfd$; by abuse of notation we denote by $m$ both
multiplications (in $A$ and $A_\mu$).
\begin{equation}\label{eq:dbbl}\xymatrix{ A_\mu A C\ar[d]_{A_\mu \psi
}\ar[r]_{A_\mu\iota C}\ar@/^1pc/[rr]^{m C}&
A_\mu A_\mu C\ar[r]_{m C}\ar[d]^{A_\mu\psi_\mu }&A_\mu C\ar[dd]^{\psi_\mu} \\
A_\mu CA\ar[r]^{A_\mu C\iota}\ar[d]_{A_\mu}&A_\mu CA_\mu\ar[d]^{\psi_\mu A_\mu}&\\
CA_\mu A\ar[r]^{CA_\mu\iota}\ar@/_1pc/[rr]_{Cm}&CA_\mu
A_\mu\ar[r]^{Cm} &CA_\mu}\end{equation}  This diagram clearly
commutes and when we tensor the whole diagram with $M$ from the
right we see that the upper left square in~(\ref{eq:doubleeq})
commutes.
$$\xymatrix{
A_\mu A CM\ar[d]_{A_\mu\psi M}\ar[r]^{A_\mu\psi M}&A_\mu
CAM\ar[r]^{A_\mu
C\nu}&A_\mu CM\ar[dd]^{\psi_\mu M}\\
A_\mu CAM\ar[d]_{\psi_\mu AM}\ar@{=}[ru]&&\\
CA_\mu A M\ar[rr]^{CA_\mu\nu}&&CA_\mu M}$$ This diagram commutes
by naturality, and shows that the lower left square
in~(\ref{eq:doubleeq}) commutes.

We conclude that the map $l$ is well-defined. We need to check
that it is a map of $A$-modules. But again, the commutativity of
the diagram~(\ref{eq:dbbl}) shows that $\psi_\mu\otimes M:
A_\mu\otimes GM\to G(A_\mu\otimes M)$ is a map of left $A$-modules
(in fact it is a map of $A_\mu$-modules simply by the pentagon for
$\psi$ and $m$); hence {\it a fortiori} the induced map on
quotients respects $A$-module structure.

\nxpoint {\bf Proposition.} {\it $l_M$ above form a distributive
law.}

{\it Proof.} Direct check.

\nxpoint {\bf Theorem.}
{\it Every localization compatible pair of entwinings
induces a continuous localization
$Q^*_\psi : {}^C_A \MM_\psi\to {}^C_{A_\mu}\MM_{\psi_\mu}$
between the categories
of entwined modules for the two entwinings
such that $U_\mu Q^*_\psi = Q^* U$ where $U$ and $U_\mu$ are
the forgetful functors from the categories of entwined to
the categories of usual modules over $A$ and $A_\mu$ respectively.
}

\section{The case of module algebras}

\nxpoint We say that the left action $\triangleright$
of a bialgebra $H$ on an algebra $A$ is
Hopf or that $(A,\triangleright)$ is a left $H$-module algebra
if $h\triangleright (a b) = \sum (h_{(1)}\triangleright a)
(h_{(2)}\triangleright a)$. Let $\iota : A\to A_\mu$ be an Ore localization
and $(A,\triangleright)$ a left $H$-module algebra. We say that $\iota$
is compatible with module algebra structure if there is a Hopf action
$\triangleright'$ of $H$ on $A_\mu$ such that 
$\iota \circ \triangleright = \triangleright'\circ (H\otimes\iota)$. 

\nxpoint {\it The monoidal category of right $H$-comodules has a canonical
action on the category of modules ${}_A\MM$ 
over a left $H$-module algebra $A$.} Again, $A$ is a monoid in that category
so we get in particular a monad ${\bf T}$ on ${}_A\MM$.

The action is induced by the distributive law with components 
$l_{M,P} : A\otimes (M\otimes P)\to (A\otimes M)\otimes P$
given by the formula $a\otimes (m\otimes p)\mapsto 
\sum (p_{(1)}\triangleright_A \otimes m)\otimes p_{(0)}$. 

One should check that one indeeds get an action. For simplicity of notation 
we do it for $P = H$; that is we check that ${\bf T}$ is a monad. General 
case is almost the same. 

Define the endofunctor $T$ on the category of left $A$-modules
\[\begin{array}{l} 
T (M,\triangleright_M) := (M \otimes_\genfd H, \triangleright_{TM}),
\,\,\,\,\,M \in A-{\rm Mod};\\
Tf := f \otimes 1_H \in {\rm Hom}_{A} (TM, TN), 
\,\,\,\,\,\forall f \in {\rm Hom}_A(M, N). 
\end{array}
\]
where the $A$-action $\triangleright_{TM}$ on $M \otimes_\genfd H$ is
given by $\genfd$-linear extension of the formula
$a \triangleright_{TM} (m\otimes h)
:= \sum  (( h_{(2)} \triangleright_A a)\triangleright_M m))\otimes h_{(1)}$.

We have to check that $\triangleright_{TM}$ is indeed an $A$-action:
\[\begin{array}{lcl}
a\triangleright_{TM} (a'\triangleright_{TM} (m\otimes h)) &=&
a\triangleright_{TM} 
(( h_{(2)} \triangleright_A a')\triangleright_M m))\otimes h_{(1)}\\
&=&  (( h_{(2)} \triangleright_A a)\triangleright_M 
((h_{(3)} \triangleright_A a')\triangleright_M m))\otimes h_{(1)}\\
\,\,(\triangleright_A \mbox{ is a Hopf action}) &= &
  ( h_{(2)} \triangleright_A (aa'))\triangleright_M m) \otimes h_{(1)}\\
&=& (aa') \triangleright_{TM} (m \otimes h).
\end{array}\]
\[\begin{array}{lcl}
1\triangleright_{TM} (m\otimes h) &=&
 (( h_{(2)} \triangleright_A 1)\triangleright_M m))\otimes h_{(1)}\\
\,\,\,\,\,\,\,\,(\triangleright_A \mbox{ is a Hopf action}) &= &
  (\epsilon( h_{(2)} ) 1\triangleright_M m))\otimes h_{(1)}\\
&=&  m \otimes h.
\end{array}\]
We also check that $Tf$ is indeed a map of left $A$-modules:
\[\begin{array}{lcl}
(Tf)[a \triangleright_{TM} (m \otimes h)] &=&
(Tf) [ (h_{(2)} \triangleright_A a) \triangleright_M m \otimes h_{(1)}]\\
&=& f((h_{(2)} \triangleright_A a) \triangleright_M m ) \otimes h_{(1)}\\
&=& (h_{(2)} \triangleright_A a) \triangleright_M f(m) \otimes h_{(1)}\\
&=& a \triangleright_{TM} (f(m)\otimes m)\\
&=& a \triangleright_{TM} [(Tf)(m \otimes h)]
\end{array}\]
Define the natural transformations $\mu : TT \Rightarrow T$ and
$\eta : {\rm Id}\Rightarrow T$ by 
\[
\mu_{(M,\triangleright_M)}(\sum_i  m_i \otimes h_i \otimes g_i)
:= \sum_i  m_i\otimes h_i g_i,
\]\[
\eta_{(M,\triangleright_M)}(m) := m\otimes 1.
\]

Here we have to check that $\mu_M := \mu_{(M,\triangleright_M)}$ and 
$\eta_M := \eta_{(M,\triangleright_M)}$ are indeed maps of left $A$-modules.
\[\begin{array}{lcl}
a \triangleright_{TTM} [(m \otimes h )\otimes g] & = &
[(g_{(2)}\triangleright_{A} a)\triangleright_{TM} (m \otimes h)]
\otimes g_{(1)}\\
&=& [h_{(2)} \triangleright_A (g_{(2)}\triangleright_A a)]
 \triangleright_m m \otimes h_{(1)}\otimes g_{(1)}\\
&=& ( h_{(2)} g_{(2)} \triangleright_A a )\triangleright_M m
 \otimes h_{(1)}\otimes g_{(1)}\\
& \stackrel{\mu_M}{\mapsto} &
( h_{(2)} g_{(2)} \triangleright_A a )\triangleright_M m
\otimes h_{(1)} g_{(1)}\\
&&=\,\,\,\, a \triangleright_{TM} (m \otimes hg)\\
&&=\,\,\,\, a \triangleright_{TM} [\mu_M ((m \otimes h)\otimes g)].
\end{array}\]
Now we have a straightforward

\nxpoint {\bf Proposition.} {\it Compatibility of Hopf action  
with localization induces a distributive law between the induced monad
${\bf T}$ defined above and the localization monad. 
}

\end{document}